\documentclass[11pt, reqno]{amsart}

\usepackage{amssymb}

\evensidemargin\oddsidemargin

\begin{document}
\baselineskip=18pt
\setcounter{page}{1}

\newcommand{\eqnsection}{
\renewcommand{\theequation}{\thesection.\arabic{equation}}
    \makeatletter
    \csname  @addtoreset\endcsname{equation}{section}
    \makeatother}
\eqnsection
   

\def\a{\alpha}
\def\B{{\bf B}}
\def\cM{{\mathcal{M}}} 
\def\cD{{\mathcal{D}}} 
\def\cG{{\mathcal{G}}} 
\def\cH{{\mathcal{H}}} 
\def\cL{{\mathcal{L}}} 
\def\cI{{\mathcal{I}}}
\def\CC{{\mathbb{C}}} 
\def\Dap{{\rm D}_{0+}^\a} 
\def\Dm{{\rm D}_{-}^\a} 
\def\Dp{{\rm D}_{+}^\a} 
\def\Ea{E_\a}
\def\esp{{\mathbb{E}}} 
\def\F{{\bf F}}
\def\Farl{{\F}_{\a,\lbd,\rho}}
\def\bF{{\bar F}}
\def\bG{{\bar G}}
\def\G{{\bf G}}
\def\Ga{{\Gamma}} 
\def\Gal{{\G_{\a,\lbd}}} 
\def\GG{{\bf \Gamma}}
\def\ii{{\rm i}} 
\def\Iap{{\rm I}_{0+}^\a} 
\def\Im{{\rm I}_{-}^\a} 
\def\Ip{{\rm I}_{+}^\a} 
\def\L{{\bf L}}
\def\lbd{\lambda}
\def\lacc{\left\{}
\def\lcr{\left[}
\def\lpa{\left(}
\def\lva{\left|}
\def\M{{\bf M}}
\def\NN{{\mathbb{N}}} 
\def\pb{{\mathbb{P}}}
\def\QQ{{\mathbb{Q}}} 
\def\R{{\bf R}}
\def\rl{{\mathbb{R}}}
\def\racc{\right\}}
\def\rpa{\right)}
\def\rcr{\right]}
\def\rva{\right|}
\def\sga{\sigma^{(\a)}}
\def\T{{\bf T}}
\def\Un{{\bf 1}}
\def\X{{\bf X}}
\def\Y{{\bf Y}}
\def\W{{\bf W}}
\def\Z{{\bf Z}}
\def\Warl{{\W}_{\a,\lbd,\rho}}
\def\Za{{\Z_\a}}

\def\elaw{\stackrel{d}{=}}
\def\claw{\stackrel{d}{\longrightarrow}}
\def\elaw{\stackrel{d}{=}}
\def\qed{\hfill$\square$}


\newtheorem{prop}{Proposition}[section]
\newtheorem{coro}[prop]{Corollary}
\newtheorem{lem}{Lemma\!\!}
\newtheorem{ex}[prop]{Example}
\newtheorem{exs}[prop]{Examples}
\newtheorem{rem}{Remark\!\!}
\newtheorem{theo}{Theorem\!\!}
\newtheorem{conj}{Conjecture\!\!}

\renewcommand{\thetheo}{}
\renewcommand{\theconj}{}
\renewcommand{\thelem}{}
\renewcommand{\therem}{}

\newcommand{\noi}{\noindent}
\newcommand{\dis}{\displaystyle }

\title{Remark on a Mittag-Leffler function of Le Roy type}

\author[T.~Simon]{Thomas Simon}

\address{Laboratoire Paul Painlev\'e, UMR 8524, Universit\'e de Lille,  Cit\'e Scientifique, F-59655 Villeneuve d'Ascq Cedex, France. {\em Email}: {\tt thomas.simon@univ-lille.fr}}

 


\begin{abstract} 

We give some necessary and some sufficient conditions for the complete monotonicity on the negative half-line of a Mittag-Leffler function of Le Roy type. It is conjectured that the underlying positive random variable, when it exists, must be logarithmically infinitely divisible. 

\end{abstract}

\maketitle

\section{Introduction}

\label{FurthA}

We consider in this letter the three-parameter entire function defined as the convergent series
$$F^{(\gamma)}_{\a,\beta} (z)\; =\; \sum_{n\ge 0} \frac{z^n}{(\Ga(\beta + \a n))^\gamma}, \qquad z\in\CC,$$
with $\a,\beta,\gamma >0.$ This function has been recently considered in \cite{G, GP, GGH, GMR, P} from various viewpoints. It is an extension of both the generalized Mittag-Leffler \cite{S} and the Le Roy function \cite{Le}, with
$$F_{\a,\beta}^{(1)}(z)\; =\; E_{\a, \beta}(z)\qquad\mbox{and}\qquad F_{1,1}^{(\gamma)}(z)\; =\; \sum_{n\ge 0} \frac{z^n}{(n!)^\gamma}\cdot$$ 
We will consider here the classical problem of the complete monotonicity on the negative half-line. It has been shown in \cite{S} that $x\mapsto E_{\a, \beta}(-x)$ is completely monotone (CM) on $(0,\infty)$ iff $\a\in [0,1]$ and $\beta \ge \a.$ It was observed in \cite{BS} - see Proposition 5.3 therein - that $x\mapsto F_{1, 1}^{(\gamma)}(-x)$ is CM on $(0,\infty)$ iff $\gamma\in [0,1]$. In the case when $\gamma$ is an integer, the problem was recently studied in \cite{GGH}, where the property is shown to be equivalent to the non-negativity on $(0,\infty)$ of a certain Meijer $G-$function. Throughout, the CM property of a given function will be implicitly meant on $(0,\infty).$ Our contribution in this note is the 

\begin{theo}
\label{MLRCM}

For every $\a,\beta,\gamma > 0,$ the following holds.

\smallskip

{\em (a)} If $\gamma \le 1,$ the function $x\mapsto F_{\a,\beta}^{(\gamma)}(-x)$ is {\em CM} iff $\a\gamma \le 1$ and $\beta \ge \a.$

\medskip

{\em (b)} If $\gamma > 1$ and $\a\gamma \neq 1,$ the function $x\mapsto F_{\a,\beta}^{(\gamma)}(-x)\;\lacc\!\!\begin{array}{l} \mbox{\!\! is {\em CM} if $\a\gamma < 1$ and $\beta \ge \a(1+\gamma)/2.$}\\
\mbox{\!\! is not {\em CM} if $\a\gamma > 1$ or $\beta \le \a.$}
\end{array}\right.$
 
\medskip

{\em (c)} If $\gamma > 1$ and $\a\gamma = 1,$ the function $x\mapsto F_{\a,\beta}^{(\gamma)}(-x)$ is {\em CM} iff $\beta \ge (1 +\a)/2.$
\end{theo}

This theorem extends in its part (a) the aforementioned characterizations of \cite{BS, S}. In its parts (b) and (c), it improves on all the results of \cite{GGH} and solves the numerical conjecture stated therein that for every $\gamma = n$ an integer, every $\a\le 1/n$ and $\beta \ge (n+1)/2n,$ the CM property holds for $x\mapsto F_{\a,\beta}^{(\gamma)}(-x).$ Part (c) shows also that the criterion on $\beta$ is a characterization when $\a = 1/n,$ whereas part (b) implies that this condition can be weakened for $\a < 1/n.$ 
  
The connection between the Stieltjes moment problem and the CM property is mentioned in the conclusion of \cite{GGH}. We follow this approach here in a continuous way, with the help of Carlson's theorem which implies an equivalence between the CM property and the existence of certain Mellin transforms. We then appeal to some classical considerations on infinite divisibility and the Gamma function. The situation $\{\gamma > 1,\, \a\gamma < 1,\, \beta \in (\a, (\a +1)/2)\}$ remains however open. A plausible, albeit non explicit, conjecture in this respect is formulated in Section 3.    

\section{Proof of the Theorem} We start with the sufficient conditions, and consider the following property:
$$\cM_{\a,\beta,\gamma}\; : \; \exists\; X\; \mbox{a positive random variable such that}\; \esp[X^s]= \Ga(1+s)\lpa\frac{\Ga(\beta)}{\Ga(\beta +\a s)}\rpa^\gamma,\; s >0.$$
It is clear that $\cM_{\a,\beta,\gamma}$ implies the CM character of $x\mapsto F_{\a,\beta}^{(\gamma)}(-x),$ since the underlying random variable has then integer moments
$$\esp[X^n]\; =\; n! \lpa\frac{\Ga(\beta)}{\Ga(\beta +\a n)}\rpa^\gamma\!\!, \;\; n\ge 0,$$
and moment generating function $\esp [e^{zX}] = (\Ga(\beta))^\gamma F_{\a,\beta}^{(\gamma)}(z),\, z\in\CC.$ In particular, by Bernstein's theorem, 
$$F_{\a,\beta}^{(\gamma)}(-x)\; =\; \frac{1}{(\Ga(\beta))^\gamma}\, \esp [e^{-xX}]\quad \mbox{is CM.}$$ 
On the other hand, the exponential formula for the Gamma function and some simplifications - see formul\ae\, 1.9(1) resp. 1.7.2(14) in \cite{EMOT} - imply
$$ \Ga(1+s)\lpa\frac{\Ga(\beta)}{\Ga(\beta +\a s)}\rpa^\gamma \; =\; \exp\lcr (\psi(1) - \a\gamma\psi(\beta))\, s\; + \; \int^{\infty}_0 (e^{-sx} - 1 + sx)\, \varphi_{\a,\beta,\gamma}(x)\, \frac{dx}{x}\rcr$$
for every $s > 0,$ where $\psi$ is the usual Digamma function and
$$\varphi_{\a,\beta,\gamma}(x)\; = \;\frac{e^{-x}}{1-e^{-x}}\, -\, \gamma\,\frac{e^{-\frac{\beta x}{\a}}}{1-e^{-\frac{x}{\a}}}$$
integrates $1\wedge x$ on $(0,\infty).$ By the L\'evy-Khintchine formula - see e.g. Chapter 1.2 in \cite{Apple} - there will hence exist a real, infinitely divisible random variable $Y$ such that 
$$\esp[e^{-sY}] \; =\;\Ga(1+s)\lpa\frac{\Ga(\beta)}{\Ga(\beta +\a s)}\rpa^\gamma\!\!, \qquad s \ge 0,$$ 
as soon as 
\begin{equation}
\label{zzz}
\varphi_{\a,\beta,\gamma}(x)\,\ge\, 0\;\;\mbox{on}\;\, (0,\infty)\quad\Longleftrightarrow\quad z \, +\, \gamma (z^{\beta -\a}  -\, z^{\beta})\, \le \, 1\;\;\mbox{for all} \;\; z\,\in\, (0,1).
\end{equation}
Putting everything together with $X= e^{-Y}
$, we see that (\ref{zzz}) implies that $F_{\a,\beta}^{(\gamma)}(-x)$ is CM.

\medskip

Suppose first $\gamma\le 1,\, \a\gamma\le 1$ and $\beta \ge \a.$ Then
$$z \, +\, \gamma (z^{\beta -\a}  -\, z^{\beta})\; \le \; z \, +\, \gamma (1  -z^{\a})\; \le \; 1$$
since the function on the right equals $\gamma$ at zero and 1 at one, and has a derivative $1-\a\gamma z^{\a -1}$ which is either positive, or negative then positive on $(0,1).$ 

Suppose next $\gamma > 1,\, \a\gamma\le 1$ and $\beta \ge \a(1+\gamma)/2.$ Then
\begin{eqnarray*}
\varphi_{\a,\beta,\gamma}(x) & = & \varphi_{1/\gamma,\frac{1+1/\gamma}{2},\gamma}(x)\; +\; \gamma\lpa \frac{e^{-\frac{(1+\gamma) x}{2}}}{1-e^{-\gamma x}}\, - \,\frac{e^{-\frac{\beta x}{\a}}}{1-e^{-\frac{x}{\a}}}\rpa \\
& \ge & \varphi_{1/\gamma,\frac{1+1/\gamma}{2},\gamma}(x)\; +\; \gamma e^{-\frac{(1+\gamma) x}{2}} \lpa \frac{1}{1-e^{-\gamma x}}\, - \,\frac{1}{1-e^{-\frac{x}{\a}}}\rpa \; \ge \;  \varphi_{1/\gamma,\frac{1+1/\gamma}{2},\gamma}(x),
\end{eqnarray*}
where the first inequality comes from $\beta \ge \a(1+\gamma)/2$ and the second one from $\a\gamma \le 1.$ Hence, by (\ref{zzz}), it remains to prove
$$z \, +\, \gamma (z^{\frac{1-1/\gamma}{2}}  -\, z^{\frac{1+1/\gamma}{2}})\, \le \, 1\;\;\mbox{for all} \;\; z\,\in\, (0,1).$$
At the maximum of the function, one has $2z + (\gamma-1) z^{\frac{1-1/\gamma}{2}}  -(\gamma +1) z^{\frac{1+1/\gamma}{2}} = 0$ and, setting $\lbd = (1+1/\gamma)/2\in (0,1),$ we are hence reduced to $z^\lbd + z^{1-\lbd}  -  z\, \le \, 1$ for all $z\in (0,1),$ which is elementary. This completes the proof of the sufficient conditions.\\ 

We now show the necessary conditions. Assuming that $x\mapsto F_{\a,\beta}^{(\gamma)}(-x)$ is CM, we first observe by analytic continuation that $z\mapsto (\Ga(\beta))^\gamma F_{\a,\beta}^{(\gamma)}(z)$ is the moment generating function on $\CC$ of a positive random variable $X$ whose positive integer moments read 
$$\esp[X^n]\; =\; n!\,\times\lpa\frac{\Ga(\beta)}{\Ga(\a n +\beta)}\rpa^\gamma,\quad n\ge 0.$$ 
If $\a\gamma > 1,$ Stirling's formula implies $\esp[X^n]^{\frac{1}{n}}\to 0$ as $n\to\infty$ so that $X\equiv 0,$ a contradiction because $F_{\a,\beta}^{(\gamma)}$ is not constant. We will assume henceforth $\a\gamma \le 1$ and we will now show that $\cM_{\a,\beta,\gamma}$ must hold. 

Assuming first $\a\gamma = 1,$ Stirling's formula implies
$$\esp[X^n]\; \sim\; \Ga(\beta)^\gamma (2\pi)^{\frac{1-\gamma}{2}} \gamma^{n + \gamma(\beta -1/2)} n^{\gamma(1/2-\beta) +1/2}\quad \mbox{as $n\to\infty.$}$$ 
In particular, the Mellin transform $s\mapsto\esp[X^s]$ is analytic on $\{\Re(s)\ge 0\},$ bounded on $\{\Re(s) = 0\},$ and has at most exponential growth on $\{\Re(s) > 0\}$ because
$$\vert \esp[X^s]\vert\; \le\; \esp\lcr X^{\Re(s)}\rcr\; =\; \lpa \esp\lcr X^{[\Re(s)] +1}\rcr\rpa^{\frac{\Re(s)}{[\Re(s)] +1}}$$
by H\"older's inequality. Since, on the other hand,  
$$\Ga(1+s)\lpa\frac{\Ga(\beta)}{\Ga(\a s +\beta)}\rpa^\gamma \; = \; \Ga(\beta)^\gamma (2\pi)^{\frac{1-\gamma}{2}} \gamma^{s + \gamma(\beta -1/2)} s^{\gamma(1/2-\beta) +1/2}\quad \mbox{as $\vert s\vert \to\infty$ with $\vert \arg s\vert \le \pi/2$},$$
reasoning exactly as in the proof of Theorem 2.3 in \cite{BS} with Carlson's theorem - see e.g. Section 5.81 in \cite{Ti} - implies the identification
\begin{equation}
\label{Ident}
\esp[X^s] \; = \; \Ga(1+s)\,\lpa\frac{\Ga(\beta)}{\Ga(\a s +\beta)}\rpa^\gamma\!, \quad s \ge 0.
\end{equation}
Observe that by analytic continuation, the identification holds also true for $s > -1.$

Assuming next $\a\gamma < 1,$ the identification (\ref{Ident}) is obtained by an analogous argument. It consists in identifying the bounded sequence
$$ \frac{1}{n!}\,\times\lpa\frac{\Ga(\a n +\beta)}{\Ga(\beta)}\rpa^\gamma$$
as the values at non-negative integer points of the function
$$ \frac{1}{\Ga(1+s)}\,\times\lpa\frac{\Ga(\a s +\beta)}{\Ga(\beta)}\rpa^\gamma\; =\; e^{-(1-\a\gamma) s (\ln (s) + O(1))}\quad\mbox{as $\vert s\vert \to \infty$ with $\vert \arg s\vert \le \pi/2.$}$$
On $\{\Re(s) \ge 0\},$ we see that this function has growth at most $e^{\pi (1-\a\gamma) \vert s\vert /2}$ and we can again apply Carlson's theorem. \\

We can now finish the proof of the necessary conditions in (a) and (b). Indeed, we see that if $\beta < \a,$ the right-hand side of (\ref{Ident}) vanishes at $s = -\beta/\a > -1,$ a contradiction since $\esp [X^{-\beta/\a}] > 0.$ For the boundary case in (b), rewriting
$$\Ga(1+s)\,\lpa\frac{\Ga(\a)}{\Ga(\a s +\a)}\rpa^\gamma\; =\; (1+s)^{\gamma -1}\, \Ga(2+s)\lpa\frac{\Ga(\a +1)}{\Ga(\a s +\a +1)}\rpa^\gamma,$$
we see that if $\gamma > 1$ and $\beta =\a,$ the random variable $X$ given by $(\Ga(\a))^\gamma F_{\a,\a}^{(\gamma)}(-x) = \esp[e^{-xX}]$ is such that $\esp [X^s] \to 0$ as $s\downarrow -1,$ another contradiction. 

The proof of the necessary condition for (c) is more involved and mimics that of Lemma 2.1 in \cite{BS}. Assuming $\{\gamma > 1, \a\gamma = 1, \beta \in (\a, (1+\a)/2)\}$ and setting $\varepsilon = (1+\a)/2 - \beta > 0,$ a Taylor expansion of order two shows first that $z  + \gamma (z^{\beta -\a}  - z^{\beta}) = 1 +\varepsilon (1-z)^2/2 + o(1-z)^2,$ whence 
\begin{equation}
\label{Neg}
\varphi_{\a,\beta,\gamma} (x) \; \rightarrow\; -\frac{\varepsilon\gamma}{2}\qquad\mbox{as $x\to 0.$}
\end{equation}
This implies that $\varphi_{\a,\beta,\gamma}$ is integrable and takes negative values on $(0,\infty).$ Moreover, the second derivative of $z\mapsto z  + \gamma (z^{\beta -\a}  - z^{\beta}) - 1$ vanishes only once, at
$$\lpa \frac{(\beta -\a)(1+\a -\beta)}{\beta (1-\beta)}\rpa^{\frac{1}{\a}}\; \in \; (0,1)$$
since $1-\beta > \beta -\a.$ By Rolle's theorem and the above Taylor expansion, this implies that the function $z\mapsto z  + \gamma (z^{\beta -\a}  - z^{\beta}) - 1$ vanishes also only once on $(0,1)$ and that it is negative then positive. It is then easy to see that there exists $x_* > 0$ such that $\varphi_{\a,\beta,\gamma} (x) < 0$ for $x\in (0, x_*)$ and $\varphi_{\a,\beta,\gamma} (x) > 0$ for $x\in (x_*, \infty).$ Now if $\cM_{\a,\beta,\gamma}$ held, since $\varphi_{\a,\beta,\gamma}$ is integrable by (\ref{Neg}), there would exist a random variable $Y_1$ with Laplace transform
$$\esp[e^{-sY_1}]\; =\; \exp\lcr\int^{\infty}_0 (e^{-sx} - 1)\, \varphi_{\a,\beta,\gamma}(x)\, \frac{dx}{x}\rcr, \qquad s > 0.$$
Considering an independent random variable $Y_2$ with Laplace transform
$$\esp[e^{-sY_2}]\; =\; \exp\lcr\int^{x_*}_0 (1-e^{-sx})\, \varphi_{\a,\beta,\gamma}(x)\, \frac{dx}{x}\rcr, \qquad s > 0,$$
it follows from Tucker's theorem \cite{Tu} that $Y_2$ is absolutely continuous since it is infinitely divisible with an absolutely continuous and, by (\ref{Neg}), infinite L\'evy measure. By e.g. Lemma 1 in \cite{Tu}, this contradicts the fact that $Y_1 + Y_2$ has an atom at zero, given by
$$\pb[Y_1 + Y_2 = 0]\; =\; \exp\lcr- \int^{\infty}_{x_*} \varphi_{\a,\beta,\gamma}(x)\, \frac{dx}{x}\rcr\; > \; 0.$$
This contradiction completes the proof of the necessary conditions.

\qed

\begin{rem}{\em Consider the following function
$$V(x)\; =\; \theta^{-x}\, \frac{(\Gamma (Ax +a))^\a}{(\Gamma (Bx +b))^\beta}$$
on $(0,\infty),$ where all the parameters are positive. With the notation of the recent paper \cite{BCK} for logarithmic completely monotone (LCM) functions, the proof of parts (a) and (c) of the theorem can be easily adapted to show that the following characterization holds:
$$V(x)\; \mbox{is CM}\quad \Longleftrightarrow\quad V(x)\; \mbox{is LCM,}$$ 
and that this property is equivalent to the following set of conditions on the parameters:
$$\a A \,=\,\beta B, \;\;\;Ab\, \ge\, Ba,\;\;\;  B^{B\beta}\theta\,\ge\, A^{A\a}\;\;\;\mbox{and}\;\;\;\beta (2b -1)\, \ge\, \a (2a - 1).$$
This also implies that the necessary conditions given in Corollary 3.3 of \cite{BCK} for the LCM property of $V(x)$ are actually sufficient in the case $p=s=1,$ with the notation therein. The latter equivalence is also an extension of Theorem 3.13 in \cite{BCK}. In particular, the function 
$$\psi(Ax +a)\, - \,\psi(Bx + b)\; =\; \sum_{n\ge 0} \lpa \frac{1}{Bx + b + n}\, -\, \frac{1}{Ax + a +n}\rpa$$ is Bernstein on $(0,\infty)$ if and only if $a\ge b$ and $2(Ab - Ba)\, \ge\, A - B \, \ge\, 0.$ We leave the details, which are available upon request, to the interested reader. See Section 3 in \cite{BCK} and the references therein for further results on the LCM property for more general ratios of powers of the Gamma function.}
\end{rem}

\section{An open problem}

The above proof shows that
$$F_{\a,\beta}^{(\gamma)}(-x)\; \mbox{is CM}\; \Longleftrightarrow\; \cM_{\a,\beta,\gamma}\; {\rm holds}\; \Longleftrightarrow\; \lacc\frac{n!}{(\Ga(\a n +\beta))^\gamma}\racc\; \mbox{is a Stieltjes moment sequence}$$
and that this is ensured by $z  + \gamma (z^{\beta -\a}  - z^{\beta}) \le 1$ for all $z\in (0,1).$ The latter criterion was also shown to be necessary for $\gamma \le 1$ or $\{\gamma > 1, \a\gamma = 1\}.$ We believe that this necessity is true in general.

\begin{conj}
\label{CMZ}
For every $\a,\beta,\gamma > 0,$ one has
$$F_{\a,\beta}^{(\gamma)}(-x)\; \mbox{is {\em CM}}\quad \Longleftrightarrow\quad z \, +\, \gamma (z^{\beta -\a}  - z^{\beta})\, \le\, 1 \quad\mbox{for all $z\in (0,1).$}$$ 
\end{conj}

In general, it does not seem that the criterion on the right-hand side can be expressed explicitly in terms of $\a,\beta,\gamma,$ save for $\gamma\le 1$ or $\{\gamma > 1, \a\gamma = 1\}$ by part (a) resp. part (c) of the theorem. If we fix $\gamma > 1$ and set
$$\beta(\a)\; =\; \inf\{ \beta > 0, \;  z \, +\, \gamma (z^{\beta -\a}  - z^{\beta})\, \le\, 1 \;\mbox{for all $z\in (0,1)$}\}\; \in\; (\a, \a(1+\gamma)/2)$$
for all $\a\in (0,1/\gamma),$ it is not difficult to prove that $\a\mapsto\beta(\a)/\a$ is non-decreasing. We also believe that $\a\mapsto \beta(\a)$ is convex, but this fact still eludes us.

As seen in Section 2, the above conjecture amounts to the property that the positive random variable $X$ associated to the CM function $F_{\a,\beta}^{(\gamma)}(-x)$ by
$$\esp[e^{zX}]\; =\; (\Ga(\beta))^\gamma F_{\a,\beta}^{(\gamma)}(z), \qquad z\in\CC,$$must be logarithmically infinitely divisible, that is the random variable $\log X$ is infinitely divisible. In the remaining situation $\{\gamma > 1, \a\gamma < 1\},$ the puzzling point is that for $\beta\in (\a, \a(1+\gamma)/2)$ the function $z\mapsto z + \gamma (z^{\beta -\a}  - z^{\beta}) - 1$ may vanish twice and be negative then positive then negative on $(0,1),$ which implies that there exists $x_2 > x_1 > 0$ such that $\varphi_{\a,\beta,\gamma} (x) < 0$ for $x\in (x_1, x_2)$ and $\varphi_{\a,\beta,\gamma} (x) > 0$ for $x\in (0, x_1)\cup (x_2, \infty).$ In this case, the absolute continuity counterargument used at the end of Section 2 does not hold anymore. Besides, there exist real random variables having characteristic function
$$z\;\mapsto\;\exp\lcr\int_{\rl} (e^{{\rm i} z x} - 1 - {\rm i} z x)\, \varphi(x)\, dx\rcr$$
for some real function $\varphi (x)$ which is negative on some open interval $(x_1,x_2)$ and non-negative otherwise. See for example Remark (c) in \cite{KSW} for a family of such random variables constructed in the framework of moments of Gamma type (which are different from those we dealt with in the present note).


\begin{thebibliography}{10}

\bibitem{Apple}
D.~Applebaum. {\em L\'evy processes and stochastic calculus.} Cambridge University Press, Cambridge, 2004.

\bibitem{BCK}
C.~Berg, A.~\c{C}etinkaya and D.~Karp. Completely monotonic ratios of basic and ordinary Gamma functions. {\tt arXiv:2004.14075} 

\bibitem{BS}
L.~Boudabsa and T.~Simon. Some properties of the Kilbas-Saigo function. {\tt arXiv:2012.05666} 

\bibitem{EMOT}
A.~Erd\'elyi, W. Magnus, F.~Oberhettinger and F.~G.~Tricomi. {\em Higher transcendental functions. Vol. I.} McGraw-Hill, New-York, 1953.

\bibitem{G}
S.~Gerhold. Asymptotics for a variant of the Mittag-Leffler function. {\em Int. Transf. Spec. Funct.} {\bf 23} (6), 397-403, 2012.

\bibitem{GP}
R.~Garra and F.~Polito. On some operators involving Hadamard derivatives. {\em Int. Transf. Spec. Funct.} {\bf 24} (10), 773-782, 2013.

\bibitem{GGH}
R.~Garrappa, K.~G\'orska and A.~Horzela. Some results on the complete monotonicity of the Mittag-Leffler functions of the Le Roy type. {\em Fract. Calc. Appl. Anal.} {\bf 22} (5), 1284-1306, 2019.

\bibitem{GMR}
R.~Garrappa, F.~Mainardi and S.~V. Rogosin. On a generalized three-parameter Wright function of the Le Roy type. {\em Fract. Calc. Appl. Anal.} {\bf 20} (5), 1196-1215, 2017.

\bibitem{KSW}
T.~Kadankova, T.~Simon and M.~Wang. On some new moments of Gamma type. {\em Stat. Probab. Letters} {\bf 165}, Article ID 108854, 2020.

\bibitem{Le}
E.~Le Roy. Valeurs asymptotiques de certaines s\'eries proc\'edant suivant les puissances enti\`eres et positives d'une variable r\'eelle. {\em Darboux Bull.} {\bf 24} (2), 245-268, 1899.

\bibitem{P}
T.~K.~Pog\'any. Integral form of Le Roy-type hypergeometric function. {\em Int. Transf. Spec. Funct.} {\bf 29} (7), 580-584, 2018.

\bibitem{S}
W.~R.~Schneider. Completely monotone generalized Mittag-Leffler functions. {\em Expo. Math.} {\bf 14}, 3-16, 1996.

\bibitem{Ti}
E.~C.~Titchmarsh. {\em The theory of functions.} Oxford University Press, Oxford, 1939.

\bibitem{Tu}
H.~G.~Tucker. Absolute continuity of infinitely divisible distributions. {\em Pacific J. Math.} {\bf 12} (3), 1125-1129, 1962.


\end{thebibliography}
\end{document}